\documentclass{elsart3-1}

\usepackage{amssymb}
\usepackage{graphicx}
\usepackage{amsmath}

\usepackage[english,francais]{babel}

\def \a {{\alpha}}
\def \b {{\beta}}
\def \d {{\delta}}
\def \l {{\lambda}}

\def \G {{\Gamma}}
\def \s {{\sigma}}

\def \R {{\mathbb {R}}}
\def \N {{\mathbb {N}}}

\def \e {{\varepsilon}}

\def \t {{\tau}}
\def \t {{\tau}}

\def \tt {{\bar t}}

\def \O {{\Omega}}
\def \phi {{\varphi}}

\def\p{\partial}

\def \è {\`e }
\def \S {\mathcal{S} }

\def \T {\mathcal{T} }

\newtheorem{theorem}{Theorem}[section]

\newtheorem{e-proposition}[theorem]{Proposition}

\newtheorem{e-definition}[theorem]{Definition\rm}

\renewcommand{\theequation}{\arabic{equation}}
\def\og{\leavevmode\raise.3ex\hbox{$\scriptscriptstyle\langle\!\langle$~}}
\def\fg{\leavevmode\raise.3ex\hbox{~$\!\scriptscriptstyle\,\rangle\!\rangle$}}

\journal{the Acad\'emie des sciences}
\begin{document}
% place in the next line the header (rubrique) chosen for your article,
% if you know it (you can also have 2, format : Header1/Header2
\centerline{Mathematical Economics}
\begin{frontmatter}


% \author{Name\thanksref{label2}}
\selectlanguage{english}
\title{Obstacle problem for Arithmetic Asian options}

% The [label1] can be suppressed if there is only one address for all authors

\selectlanguage{english}
\author[authorlabel2]{Laura Monti},
\ead{monti@unibo.it}
\author[authorlabel2]{Andrea Pascucci}
\ead{pascucci@dm.unibo.it}

\address[authorlabel2]{Dipartimento di Matematica, Universit\`{a}
di Bologna, Piazza di Porta S. Donato 5, 40126 Bologna (Italy)}

% If you know the dates of reception, and acceptation you can put them now;
%  idem the name of the person presenting the Note

\medskip
\begin{center}
{\small Received *****; accepted after revision +++++\\
Presented by £££££}
\end{center}

\begin{abstract}
\selectlanguage{english}
% Text of abstract in English
We prove existence, regularity and a Feynman-Ka\v{c} representation formula of the strong solution
to the free boundary problem arising in the financial problem of the pricing of the American Asian
option with arithmetic average.

%{\it To cite this article: A. Name1, A.
%Name2, C. R. Acad. Sci. Paris, Ser. I 340 (2005).}

\vskip 0.5\baselineskip

\selectlanguage{francais}
% Text of abstract in French
\noindent{\bf R\'esum\'e} \vskip 0.5\baselineskip \noindent {\bf Probl\`eme de l'obstacle pour
l'option am\'ericain asiatique \`a moyenne arithm\'etique.} On d\'emontre l'existence, la
r\'egularit\'e et une formule de repr\'esentation de Feynman-Ka\v{c} de la solution forte d'un
probl\`eme avec fronti\`ere libre. Ce type de probl\`eme on le retrouve en finance pour \'evaluer
le prix d'une option asiatique \`a moyenne arithm\'etique de style am\'ericain.

%On d\'emontre l'existence, la r\'egularit\'e et une formule de repr\'esentation de Feynman-Ka\v{c}
%de la solution forte du probl\'eme \`a fronti\'ere libre qui se pr\'esente au sein du probl\'eme
%financi\'ere de la valutation de option asiatique \`a moyenne arithm\'etique de style am\'ericain.
%
%
%
%Nous alons prouver l'existence, la r\'egularit\'e et une formule de repr\'esentation de
%Feynman-Ka\v {c} pour la solution forte d'un probl\'eme avec fronti\'ere libre. Ce type de
%probl\'eme on le retrouve en finance pour \'evaluer le prix d'une moyenne arithm\'etique d'option
%asiatique de type am\'ericain.

%{\it Pour citer cet article~: A. Name1, A. Name2, C. R. Acad. Sci.
%Paris, Ser. I 340 (2005).}

\end{abstract}
\end{frontmatter}

% now the Version française abrégée, if it exists
%\selectlanguage{francais}
%\section*{Version fran\c{c}aise abr\'eg\'ee}
% Text of your Version française abrégée here.
% Note you do not need to repeat here equations that you use in the
% main text - for example 'voir (3)' is quite acceptable.

\selectlanguage{english}
% main text
\section{Introduction}%\label{}

According to the classical financial theory (see, for instance, \cite{Peskir06}) the study of
Asian options of American style leads to
%optimal stopping problems equivalent to
free boundary problems for degenerate parabolic PDEs. More precisely, let us assume that, in the
standard setting of local volatility models, the dynamics of the underlying asset is driven by the
SDE
\begin{equation}\label{1}
  dS_{t}=\mu(t,S_{t}) S_{t}dt+\s(t,S_{t})S_{t}dW_{t},
\end{equation}
and consider the process $dA_{t}=f(S_{t})dt$, where $f(S)=S$ and $f(S)=\log S$ occur respectively
in the study of the Arithmetic average and Geometric average Asian options. Then the price of the
related Amerasian option with payoff function $\phi$ is the solution of the obstacle problem with
final condition
\begin{align}\label{2}
 \begin{cases}
  \max\{L  u,\phi-u\}=0,\qquad  & ]0,T[\times \R^{2}_{+},\\
  u(T,s,a)=\phi(T,s,a), & s,a>0,
  \end{cases}
\end{align}
where
\begin{equation}\label{3}
  L u=\frac{\s^{2}s^{2}}{2}\, \p_{ss}u+r s\p_{s}u+f(s)\p_{a}u+\p_{t}u-ru,\qquad s,a>0,
\end{equation}
is the Kolmogorov operator of $(S_{t},A_{t})$ and $r$ is the risk free rate.
%We also recall that problems in the form \eqref{2}
Recently this problem has also been considered in the study of pension plans in \cite{Friedman02}
and stock loans in \cite{Dai09}.

Typical Arithmetic average payoffs are of the form
\begin{equation}\label{11}
\begin{split}
  \phi(t,s,a)=&\left(\frac{a}{t}-K\right)^{+}\qquad \text{(fixed strike)},\\
  \phi(t,s,a)=&\left(\frac{a}{t}-s\right)^{+}\qquad \text{(floating strike)}.
\end{split}
\end{equation}
A direct computation shows that in these cases a super-solution\!\!\footnote{$\bar{u}$ is a
super-solution of \eqref{2} if $$\begin{cases}
  \max\{L \bar{u},\phi-\bar{u}\}\le 0,\qquad  & ]0,T[\times \R^{2}_{+},\\
  \bar{u}(T,s,a)\ge\phi(T,s,a), & s,a>0.
 \end{cases}$$}
to \eqref{2} with $f(s)=s$ is given by
\begin{equation}\label{5}
  \bar{u}(t,s,a)=\frac{\a}{t}\left(1+e^{-\b t}\sqrt{s^{2}+a^{2}}\right)
\end{equation}
for $\a,\b$ are positive constants, with $\b$ suitably large. On the other hand it is well-known
that generally \eqref{2} does not admit a smooth solution in the classical sense.

Recently the Geometric Asian option has been studied under the following hypotheses:

  {(H1)} $\s$ is bounded, locally H\"older continuous and such that $\s\ge\s_{0}$ for some positive constant
  $\s_{0}$;

  {(H2)} $\phi$ is %continuous on $]0,T]\times \R^{2}_{+}$,
  locally Lipschitz continuous on $]0,T[\times \R^{2}_{+}$ and the distributional derivative $\p_{ss}\phi$ is
  locally

  \qquad\ lower bounded (to fix ideas, this includes $\phi(s)=(s-K)^{+}$ and excludes $\phi(s)=-(s-K)^{+}$).

\noindent In \cite{DifPa2008}, \cite{Pascucci2008} it is proved that problem \eqref{2}, with
$f(s)=\log s$, has a strong solution $u$ in the Sobolev space
\begin{equation}\label{4}
  \S^{p}_{\text{\rm loc}}=\{u\in L^{p}\mid \p_{s}u,\p_{ss}u,(f(s) \p_{a}+\p_{t})u\in L^{p}_{\text{loc}}\},
  \qquad p\ge 1.
\end{equation}
%and is continuous on $]0,T]\times \R^{2}_{+}$.
Moreover uniqueness has been proved via Feynman-Ka\v{c} representation. However, as we shall see
below, the Geometric and Arithmetic cases are structurally quite different.

The aim of this note is to give an outlined proof of the following
\begin{theorem}\label{t1}
Consider problem \eqref{2} with $f(s)=s$, under the assumptions (H1) and (H2). Then we have:

i) if there exists a super-solution $\bar{u}$ then there also exists a strong solution $u\in
\S^{p}_{\text{\rm loc}}\cap C\left(]0,T]\times \R^{2}_{+}\right)$ for any $p\ge 1$, such that $u
\le \bar{u}$;

ii) if $u$ is a strong solution to \eqref{2} such that
\begin{equation}\label{6}
  |u(t,s,a)|\le \frac{C}{t} (1+s^{q}+a^{q}),\qquad s,a>0,\ t\in]0,T],
\end{equation}
for some  positive constants $C,q$, then
\begin{equation}\label{7}
  u(t,s,a)=\sup\limits_{\t\in\T_{t,T}}E\left[e^{-r \t}\phi(\t,S^{t,s,a}_{\t},A^{t,s,a}_{\t})\right],\qquad
  t,s,a>0,
\end{equation}
where $\T_{t,T}=\{\t\in\T\mid \t\in[t,T] \text{ a.s.}\}$ and $\T$ is the set of all stopping times
with respect to the Brownian filtration. In particular there exists at most one strong solution of
\eqref{2} verifying the growth condition \eqref{6}.

%or every compact subset $M\subset\subset]0,T[\times \R^{2}_{+}$ there exists a real (not
%necessarily positive) constant $C$ such that $$
\end{theorem}

\smallskip For simplicity here we only consider the case of constant $\s$. Then by a transformation (cf.
formula (4.4) in \cite{BPV2001}), operator $L$ in \eqref{3} with $f(s)=s$ can be reduced in the
canonical form
\begin{equation}\label{8}
  L_{A}=x_{1}^{2}\p_{x_{1}x_{1}}+x_{1}\p_{x_{2}}+\p_{t},\qquad x=(x_{1},x_{2})\in\R^{2}_{+}.
\end{equation}
Before proceeding with the proof, we make some preliminary comments.

\smallskip\noindent{\bf Remark 1.} In the Geometric case $f(s)=\log s$, by a logarithmical change of variables the pricing
operator $L$ takes the form
\begin{equation}\label{9}
  L_{G}=\p_{x_{1}x_{1}}+x_{1}\p_{x_{2}}+\p_{t},\qquad x=(x_{1},x_{2})\in\R^{2}.
\end{equation}
It is known (cf. \cite{LanconelliPolidoro}) that $L_{G}$ has remarkable invariance properties with
respect to a homogeneous Lie group structure: precisely, $L_{G}$ is invariant with respect to the
left translation in the law $(t',x')\circ(t,x)=(t'+t, x'_{1}+x_{1},x'_{2}+x_{2}+t x'_{1})$ and
homogeneous of degree two with respect to the dilations $\d_{\l}(t,x)=(\l^{2}t,\l
x_{1},\l^{3}x_{2})$ in the sense that, setting $z=(t,x)$,
  $$L_{G}(u(z'\circ z))=(L_{G}u)(z'\circ z),\qquad L_{G}(u(\d_{\l}(z)))=\l^{2}(L_{G}u)(\d_{\l}(z)),\qquad z,z'\in\R^{3},\ \l>0.$$
Moreover $L_{G}$ has a fundamental solution $\G_{G}(t,x;T,X)$ (here $(t,x)$ and $(T,X)$ represent
respectively the starting and ending points of the underlying stochastic process) of Gaussian type
whose explicit expression is known explicitly:
  $$\G_{G}(t,x;T,X)=\G_{G}\left((T,X)^{-1}\circ(t,x);0,0\right)$$
where $(T,X)^{-1}=(-T,-X_{1},-X_{2}+TX_{1})$ and
%  $$\G_{G}\left(t,x,y;0,0,0\right)=\frac{\sqrt{3}}{2\pi t^{2}}\exp\left(\frac{x^{2}}{t}+\frac{3x(y-t x)}{t^{2}}
%  +\frac{3(y-t x)^{2}}{t^{3}}\right),\qquad t<0,\ x,y\in\R.$$
  $$\G_{G}\left(t,x;0,0\right)=\frac{\sqrt{3}}{2\pi t^{2}}\exp\left(\frac{x_{1}^{2}}{t}+\frac{3x_{1}(x_{2}-t x_{1})}{t^{2}}
  +\frac{3(x_{2}-t x_{1})^{2}}{t^{3}}\right),\qquad t<0,\ x_{1},x_{2}\in\R.$$
On the contrary, the Arithmetic operator $L_{A}$ does not admit a homogeneous structure:
nevertheless in this case we are able to find an interesting invariance property with respect to
the ``translation'' operator
%\begin{equation}\label{10}%(t_{0},x_{0},t_{0})\ast(t,x,y)
%  \ell_{(t_{0},x_{0},y_{0})}(t,x,y)=\left(t_{0}+t,x_{0}x,y_{0}+x_{0}y\right), \qquad t_{0},t\in\R,\
%  x_{0},y_{0},x,y>0;
%\end{equation}
\begin{equation}\label{10}%(t_{0},x_{0},t_{0})\ast(t,x,y)
  \ell_{(t',x')}(t,x)=\left(t'+t,x'_{1}x_{1},x'_{2}+x'_{1}x_{2}\right), \qquad t,t'\in\R,\
  x,x'\in\R^{2}_{+};
\end{equation}
precisely we have
\begin{equation}\label{12}%L_{A}(u(z_{0}\ast z))=\left(L_{A}u\right)(z_{0}\ast z)
 L_{A}(u\left(\ell_{(t',x')}\right))=(L_{A}u)\left(\ell_{(t',x')}\right).
\end{equation}
We remark that the fundamental $\G_{A}$ of $L_{A}$ is not known explicitly and one of the key
point in the proof consists in showing suitable summability properties of $\G_{A}$ near the pole.
Note also that, due to the invariance property \eqref{12}, the fundamental solution $\G_{A}$
%(i.e. the transition density of the process $(S_{t},A_{t})$)
verifies
  $$\G_{A}(t,x;T,X)=\frac{1}{X_{1}^{2}}\G_{A}\left(t-T,\frac{x_{1}}{X_{1}},\frac{x_{2}-X_{2}}{X_{1}};0,1,0\right),\qquad x_{1},x_{2},X_{1},X_{2}>0,\ x_{2}>X_{2},\ t<T.$$
%  $$\G_{A}(t,x,y;T,X,Y)=\frac{1}{X^{2}}\G_{A}\left(t-T,\frac{x}{X},\frac{y-Y}{X};0,1,0\right),\qquad x,y,X,Y>0,\ y>Y,\ t<T.$$

\smallskip\noindent{\bf Remark 2.} Given a domain $\O$ and $\a\in]0,1[$, we denote by $C^{\a}_{G}(\O)$ and $C^{1,\a}_{G}(\O)$
the intrinsic H\"older spaces defined by the norms
  $$\|u\|_{C^{\a}_G(\O)}=\sup_{\O}|u|+ \sup_{\stackrel{z,z_{0} \in \O}{z \neq z_{0}}}
    \frac{|u(z) - u(z_{0})|}{\|z_{0}^{-1}\circ z\|_{G}^{\a}},\qquad \|u\|_{C^{1,\a}_G(\O)}= \|u\|_{C^{\a}_G(\O)} + \|\p_{x_{1}}u\|_{C^{\a}_G(\O)}$$
where $\|\cdot\|_{G}$ is a $\d_{\l}$-homogeneous norm in $\R^{3}$. Due to the embedding Theorem
2.1 in \cite{DifPa2008}, the strong solution $u$ of Theorem \ref{t1} is locally in $C^{1,\a}_{G}$
for any $\a\in]0,1[$. In particular this ensures the validity of the standard ``smooth pasting''
condition in the $x_{1}$ variable (corresponding to the asset price). We also recall that the
(optimal) $\S^{\infty}_{\text{loc}}$-regularity of the solution has been recently proved in
\cite{FNPP}.

\smallskip\noindent{\bf Remark 3.} In the standard Black \& Scholes setting and for particular homogeneous payoffs,
it is possible to reduce the {\it spatial} dimension of the pricing problem (from two to one). In
the case of European Asian options, this was first suggested in \cite{Ingersoll} and
\cite{RogersShi} respectively for the floating and fixed strike payoffs. In the American case, the
dimensional reduction is only possible for the floating strike payoff while the fixed strike
Amerasian option necessarily involves a 2-dimensional degenerate PDE.

\smallskip\noindent{\bf Proof of Theorem \ref{t1}.} It is not restrictive to assume that $\phi$ is continuous on $[0,T]\times \R^{2}_{+}$
or equivalently we may study the problem on $[\e,T]\times \R^{2}_{+}$ for a fixed, but arbitrary,
positive $\e$. We divide the proof in some steps.

\smallskip\noindent{\bf Step 1.} Let $D_{r}(x)$ denote the Euclidean ball centered at $x\in\R^{2}$, with radius
$r$. We construct a sequence of ``lentil shaped'' domains
$O_{n}=D_{n}\left(n+\frac{1}{n},0\right)\cap D_{n}\left(0,n+\frac{1}{n}\right)$ covering
$\R^{2}_{+}$. For any $n\in \N$, the cylinder $H_{n}=]0,T[\times O_{n}$ is a {\it $L_{A}$-regular
domain} in the sense that there exists a barrier function at any point of the parabolic boundary
$\p_{P}H_{n}:=\p H_{n}\setminus(\{0\}\times O_{n})$ and therefore the obstacle problem on $H_{n}$ %\cap\{t>\e\}
has a strong solution % for any $\e>0$
(cf. Theorem 3.1 in \cite{DifPa2008}). Indeed on any compact subset $H$ of $\R\times\R^{2}_{+}$,
we have $L_{A}\equiv L_{H}$ where
\begin{equation}\label{15}
 L_{H}=a_{H}(t,x)\p_{x_{1}x_{1}}+x_{1}\p_{x_{2}}+\p_{t}
\end{equation}
and $a_{H}$ is a some smooth function such that $0<\underline{a}_{\,H}\le a_{H}\le \bar{a}_{H}$ on
$]0,T[\times\R^{2}_{+}$, with $\underline{a}_{\,H},\bar{a}_{H}$ suitable positive constants. Note
that $L_{H}$ is a perturbation of the Geometric operator $L_{G}$: as such, by Theorem 1.4 in
\cite{AMRX}, it has a fundamental solution $\G_{H}$ that is bounded from above and below by
Gaussian functions. Then by Theorem 3.1 in \cite{DifPa2008} we have: for any $n\in\N$ and $g\in
C(H_{n}\cup\p_{P}H_{n})$, $g\ge \phi$, there exists a strong solution\footnote{Here
$\S^{p}_{\text{\rm loc}}$ is the Sobolev space defined in \eqref{4} with $f(s)=s$ and $s=x_{1}$,
$a=x_{2}$.} $u_{n}\in \S^{p}_{\text{\rm loc}}\left(H_{n}\right)\cap
C\left(H_{n}\cup\p_{P}H_{n}\right)$ to problem
\begin{equation}\label{e1}
  \begin{cases}
   \max\{L_{A} u,\phi-u\}=0\  &\text{in } H_{n},\\
   u\vert_{\p_{P}H_{n}}=g. &
   \end{cases}
\end{equation}
Moreover, for every $p\ge 1$ and $H\subset\subset H_{n}$ there exists a positive constant $C$,
only dependent on $H,H_{n},p,\|\phi\|_{L^{\infty}(H_{n})},\|g\|_{L^{\infty}(H_{n})}$ such that
$\|u_{n}\|_{\S^{p}(H)}\le C$.

\smallskip\noindent{\bf Step 2.} To prove part {\it i)}, we consider a sequence of
cut-off functions $\chi_{n}\in C_{0}^{\infty}(\R^{2}_{+})$ such that $\chi_{n}=1$ on $O_{n-1}$,
$\chi_{n}=0$ on $\R^{2}_{+}\setminus O_{n}$ and $0\le \chi_{n}\le 1$. We set
$g_{n}(t,x,y)=\chi_{n}(x,y)\phi(t,x,y)+(1-\chi_{n}(x,y))\bar{u}(t,x,y)$ and denote by $u_{n}$ the
strong solution to \eqref{e1} with $g=g_{n}$. By the comparison principle we have $\phi\le
u_{n}\le u_{n+1}\le \bar{u}$ and therefore, by the a priori estimate in $\S^{p}$, for every $p\ge
1$ and $H\subset\subset H_{n}$ we have $\|u_{n}\|_{\S^{p}(H)}\le C$ for some constant $C$
dependent on $H$ but not on $n$. Then we can pass to the limit as $n\to \infty$, on compacts of
$]0,T[\times\R_{+}^{2}$, to get a solution of $\max\{L  u,\phi-u\}=0$. A standard argument based
on barrier functions shows that $u(t,x)$ is continuous up to $t=T$ and attains the final datum.
This concludes the proof of part {\it i)}.

\smallskip\noindent{\bf Step 3.} To prove part {\it ii)}, we first construct the fundamental solution $\G_{A}$ of $L_{A}$ as the
limit of an increasing sequence of Green functions for $L_{A}$ on the cylinders $H_{n}$, $n\in\N$:
to this end we combine some classical PDE technique (cf. Chapter 15 in \cite{Friedman76}) with the
recent interior and boundary Schauder estimates for degenerate Kolmogorov operators proved in
\cite{NPP}. We also show that $\G_{A}$ is the transition density of the underlying stochastic
process.

Next we observe that, for fixed $(t,x)\in ]0,T[\times \R^{2}_{+}$ and $n\in\N$ suitably large so
that $x\in O_{n}$, by the maximum principle we have
  $$\G_{A}(t,x;\cdot,\cdot)\le \G_{H_{n}}(t,x;\cdot,\cdot)+\max_{[t,T]\times \p O_{n}}\G_{A}(t,x;\cdot,\cdot)\qquad \text{ in }]t,T]\times O_{n}.$$
Consequently we infer (cf. formula (4.8) in \cite{Pascucci2008}) that $\G_{A}(t,x;\cdot,\cdot)\in
L^{p}(H_{n})$ for some $p>1$. This local summability property of $\G_{A}$ can be combined with the
standard maximal estimate
  $$E\left[\sup_{t\le \t\le T}|X^{t,x}_{\t}|^{q}\right]<\infty, \qquad q\ge 1,$$
valid for the solution $X$ of a SDE whose coefficients have at most linear growth: then, adapting
the arguments used in the proof of Theorem 4.3 in \cite{Pascucci2008}, we can show formula
\eqref{7}. This concludes the proof of part {\it ii)}.\hfill$\Box$

% \appendix

% \section{}
% \label{}

% The Acknowledgements are an un-numbered section
%\section*{Acknowledgements}
% Acknowledgements text here

%\begin{thebibliography}{00}

\bibliography{bib}

%% please try to use the bibitem system -
%% the references should be in alphabetical order of authors' names.
%% Articles with a single author first, author will 1 co-author next,
%% then author with several co-authors;
%\bibitem{CorielliPascucci} Corielli, F. and Pascucci, A. {\it Parametrix approximations for diffusion transition densities},
%preprint available at \texttt{http://www.dm.unibo.it/matecofin/member.php?id=6\&m=pascucci}
%
%\bibitem{CoxRoss1976} Cox, J. and Ross, S. {\it The valuation of options for alternative stochastic processes},
%   J. Financial Economics, 3, (1976) pp.{145--166}
%
%\bibitem{FosPierPol2008} Foschi, P. Pieressa, L. and Polidoro, S. {\it Parametrix approximations for non constant coefficient parabolic {PDE}s},
%MPRA Working Paper n.7943 (2008)
%
%\bibitem{Friedman} Friedman, A. {\it Partial differential equations of parabolic type},
%Prentice-Hall Inc., Englewood Cliffs, N.J., (1964)
%
%
%\bibitem{Levi1907} Levi, E. E., {\it Sulle equazioni lineari totalmente ellittiche alle derivate parziali},
%Rend. Circ. Mat. Palermo 24, (1907) pp.{275--317}
%

%\end{thebibliography}

\end{document}